\documentclass[12pt]{article}

\usepackage{amsmath, graphicx, color, amsthm, amsfonts, amssymb, tikz}
\usepackage{lipsum}
\usepackage{mwe}
\usepackage{lineno}

\def \bsk {\bigskip}

\newtheorem{Theorem}{Theorem}[section]

\newtheorem{Question}[Theorem]{Question}
\newtheorem{Definition}[Theorem]{Definition}
\def \bdf {\begin{Definition}\rm }
\def \edf {\end{Definition}}
\def \btm {\begin{Theorem}}
\def \etm {\end{Theorem}}
\newtheorem{Lemma}[Theorem]{Lemma}
\def \blm {\begin{Lemma}}
\def \elm {\end{Lemma}}
\newtheorem{Proposition}[Theorem]{Proposition}
\def \bpn {\begin{Proposition}}
\def \epn {\end{Proposition}}
\newtheorem{Corollary}[Theorem]{Corollary}

\def \bcr {\begin{Corollary}}
\def \ecr {\end{Corollary}}
\newtheorem{Example}[Theorem]{Example}
\def \bex {\begin{Example}\rm }
\def \eex {\end{Example}}
\newtheorem{Remark}[Theorem]{Remark}
\def \brm {\begin{Remark}\rm }
\def \erm {\end{Remark}}
\newtheorem{Observation}[Theorem]{Observation}
\def \bob {\begin{Observation}\rm }
\def \eob {\end{Observation}}

\newtheorem{Conjecture}[Theorem]{Conjecture}
\def \bc {\begin{Conjecture}}
\def \ec {\end{Conjecture}}
\newtheorem{Problem}[Theorem]{Problem}
\def \bpm {\begin{Problem}}
\def \epm {\end{Problem}}
\def \hsp {\hspace{0.5em}}
\def \bpf {\begin{proof}}
\def \epf {\end{proof}}
\def \hsp {\hspace{0.5em}}

\def \soc {\chi_\mathrm{so}}

\def \vp {\varphi}
\def \empt {\varnothing}

\def\lf{\left\lfloor}
\def\rf{\right\rfloor}

\begin{document}

\title{On strong odd colorings of graphs}
\author{Yair Caro\thanks{\hsp Department of Mathematics,
  University of Haifa-Oranim, Tivon 36006, Israel}\,,
 Mirko Petru\v sevski\thanks{\hsp Faculty of Mechanical Engineering, Ss. Cyril and Methodius University in Skopje, Macedonia}\,,
  Riste \v Skrekovski\thanks{\hsp Faculty of Mathematics and Physics, University of Ljubljana \& Faculty of Information Studies, Novo mesto, Slovenia}\,,
   Zsolt Tuza\thanks{\hsp HUN-REN Alfr\'ed R\'enyi Institute of Mathematics,
    H–1053 Budapest, Re\'altanoda u.~13--15, Hungary;
     and
    Department of Computer Science and Systems Technology,
    University of Pannonia, 8200 Veszpr\'em, Egyetem u.~10, Hungary}}
    \date{ }
\maketitle


\begin{abstract}
A strong odd coloring of a simple graph $G$ is a proper coloring of the vertices of $G$ such that for every vertex $v$ and every color $c$, either $c$ is used an odd number of times in the open neighborhood $N_G(v)$ or no neighbor of $v$ is colored by $c$. The smallest integer $k$ for which $G$ admits a strong odd coloring with $k$ colors is the strong odd chromatic number, $\soc(G)$. These coloring notion and graph parameter were recently defined in [H.~Kwon and B.~Park, \textit{Strong odd coloring of sparse graphs}, ArXiv:2401.11653v2]. We answer a question raised by the originators concerning the existence of a constant bound for the strong odd chromatic number of all planar graphs. We also consider strong odd colorings of trees, unicyclic graphs and graph products.
\end{abstract}

\medskip



\newpage

\section{Introduction}

We use standard terminology according to~\cite{BonMur08} except for few notations defined in the text. 
All considered graphs are finite, loopless and undirected. If there are no parallel edges we speak of a \textit{simple} graph. Let $G=(V_G,E_G)$ be a simple graph with vertex set $V_G$ and edge set $E_G$. For a vertex $v$, we use $N_G(v)$ to denote the open neighborhood of $v$, and $\deg_G(v)$ to denote the degree of $v$, i.e., the number of edges of $G$ incident with $v$. A \textit{$k$-vertex}, a \textit{$k^-$-vertex}, and a \textit{$k^+$-vertex}, respectively, is a vertex of degree $k$, of degree at most $k$, and of degree at least $k$. Let $\Delta(G)$ denote the maximum degree of $G$. The maximum average degree $\mathrm{mad}(G)$ of $G$ is the maximum of $\frac{2|E_H|}{|V_H|}$ over all non-void subgraphs $H$ of $G$. The girth, $g(G)$, is the length of a shortest cycle in $G$. 
A graph (as a combinatorial object) is said to be \textit{planar} if it can be
embedded in the Euclidean plane (or equivalently on the sphere) without edges crossing; such a drawing is called a \textit{plane} graph (a geometrical object).

A $k$-(vertex-)coloring of a graph $G$ is an assignment $\varphi: V_G\to\{1,\ldots,k\}$. A coloring $\varphi$ is said to be \textit{proper} if every color class (i.e., preimage $\varphi^{-1}(c)$ of any color $c$) is an independent subset of the vertex set of $G$.
 The minimum $k$ for which a graph $G$ admits a proper $k$-coloring is the \textit{chromatic number} of $G$, denoted $\chi(G)$. A \textit{square coloring} of a graph $G$ is a proper coloring of $G^2$, the latter being the graph obtained from $G$ by linking with an edge each pair of non-adjacent vertices that have a common neighbor (in other words, vertices that are at distance $2$ in $G$). 
 In 1977, Wegner conjectured an upper bound for $\chi(G^2)$ when $G$ is planar.

\begin{Conjecture}[Wegner~\cite{Weg77}]
    \label{conj:Weg}
Let $G$ be a planar graph. Then
\begin{equation*}
\chi(G^2)\leq
\begin{cases}
7 & \text{\quad if\, } \Delta(G)=3\,,\\
\Delta(G)+5 & \text{\quad if \,} 4\leq \Delta(G)\leq7\,,\\
\lf\frac{3\Delta(G)}{2}\rf+1 & \text{\quad if \,} \Delta(G)\geq8.
\end{cases}
\end{equation*}
\end{Conjecture}

Thomassen~\cite{Tho18} and Hartke et al.~\cite{HarJahTho16} independently proved Conjecture~\ref{conj:Weg} when $\Delta(G)=3$. For related research on square coloring problems and best known bounds on $\chi(G^2)$ see \cite{BorIva12, BouDesMeyPie23, BuSha16, BuZhu12, Cra22, Den22, Den23, DonXu17, DvoSkrTan08, La22}.

A recently introduced intermediate coloring concept, which can be seen as both a weakening of square coloring and a strengthening of proper coloring, has attracted considerable attention among researchers. Defined by Petru\v{s}evski and \v{S}krekovski~\cite{PetSkr22}, an \textit{odd coloring} of $G$ is a proper coloring such that for any non-isolated vertex $v$ in $G$, there is a color that is used
an odd number of times on the neighborhood of $v$. The \textit{odd chromatic number} of $G$, denoted $\chi_o(G)$, is the minimum number of colors in an odd coloring of $G$. Determining the odd chromatic number of a graph
and deciding whether a bipartite graph is odd $k$-colorable for $k\geq3$ are two $\mathcal{NP}$-complete problems~\cite{AhnImOum22}. Petru\v{s}evski and \v{S}krekovski~\cite{PetSkr22} showed
that all planar graphs are odd 9-colorable, and conjectured the following.

\begin{Conjecture}[Petru\v{s}evski and \v{S}krekovski~\cite{PetSkr22}]
    \label{conj:PetSkr}
For every planar graph $G$, it holds that $\chi_o(G)\leq5$.
\end{Conjecture}

Note that a $5$-cycle is a planar (moreover, an outerplanar) graph whose odd chromatic number is exactly $5$.  A study on various
aspects of the odd chromatic number was carried out by Caro, Petru\v{s}evski and \v{S}krekovski~\cite{CarPetSkr22}. Building upon some results in~\cite{CarPetSkr22}, Petr and Portier~\cite{PetPor23} proved that planar graphs are odd 8-colorable. 
Besides planar graphs with bounded girth, some special embedded graphs have also been investigated, such as 1-planar graphs~\cite{CraLafSon23,LiuWanYu23}, outerplanar graphs~\cite{CarPetSkr22}, outer-1-planar graphs and 2-boundary planar graphs~\cite{QiZha22}, and $k$-planar graphs~\cite{Hic22,DuiMorOda22}.  
Tian and Yin~\cite{TiaYin22} showed that every toroidal graph with girth at least 4 is odd 7-colourable. Odd colorings of sparse graphs, measured in terms of the maximum average degree, were investigated by Cranston~\cite{Cran22} and Cho, Choi, Kwon, and Park~\cite{ChoChoKwoPar232}. 
Their findings include that
every planar graph of girth at least $5$, $7$ or $11$, respectively, is odd $6$-, $5$- or $4$-colorable. 
Most recently, Miao, Sun, Tu and Yu~\cite{MiaSunTuYu24} gave support to Conjecture~\ref{conj:PetSkr} by showing that all triangle-free planar graphs without intersecting $4$-cycles are odd $5$-colorable. They also showed that planar graphs with non-intersecting (resp.\ edge-disjoint) 3-cycles are odd $6$-colorable (resp.\ odd $7$-colorable).

\medskip

In this paper, we delve into a natural stronger variant of the odd coloring notion, introduced of late by  imposing a stronger neighborhood constraint on the properness condition. Defined by Kwon and Park~\cite{KwoPar24}, a \textit{strong odd coloring} of a graph $G$ is a proper coloring of its vertices such that for every non-isolated vertex $v$, if a color $c$ is used on $N_G(v)$ then $c$ appears an odd number of times on the vertices comprising $N_G(v)$. Clearly, the proper setting for any study of this coloring concept are the simple graphs, since parallel edges are irrelevant for either requirement. Every simple graph admits a strong odd coloring since a square coloring is also a strong odd coloring. The \textit{strong odd chromatic number} of $G$, denoted $\soc(G)$, is the minimum $k$ such that $G$ admits a strong odd $k$-coloring. It is clear from the definitions that for every graph $G$,

\begin{equation*}
\chi(G)\leq\chi_o(G)\leq\soc(G)\leq\chi(G^2)\leq(\Delta(G))^2+1\,.
\end{equation*}

Clearly, for every claw-free graph $G$ it holds that $\soc(G)=\chi(G^2)$. However, as observed in~\cite{KwoPar24}, the difference $\chi(G^2)-\soc(G)$ can be arbitrarily large: indeed, if $G=K_{m,n}$ then $\chi(G^2)=m+n$ whereas $\soc(G)\leq4$. The originators of the concept investigated strong odd colorings of sparse graphs $G$. They proved that if $G$ has $\mathrm{mad}(G)\leq\frac{20}{7}$ then $\soc(G)\leq\Delta(G)+4$; moreover, the bound is tight for subcubic planar graphs. Also if $G$ is a graph with $\mathrm{mad}(G)\leq\frac{30}{11}$ and $\Delta(G)\geq4$ then $\soc(G)\leq\Delta(G)+3$. 
Since results regarding bounded maximum average degree have natural corollaries to planar graphs
with girth restrictions, because the class of graphs $G$ with $\mathrm{mad}(G)<\frac{2g}{g-2}$ includes all planar graphs with girth at least $g$, the findings of this pair of scholars yield that every planar graph $G$ with girth $g\geq7$ (resp.\ $g\geq8$) has $\soc(G)\leq\Delta(G)+4$ (resp.\ $\soc(G)\leq\Delta(G)+3$).
In view of the odd 8-colorability of all planar graphs~\cite{PetPor23},  Kwon and Park ended~\cite{KwoPar24} with the following natural question.

\begin{Question}
    \label{q:constant}
Is there a constant $C$ such that every planar graph $G$ has $\soc(G)\leq C$?
\end{Question}

One of the main contributions of the present work is a proof that such a finite upper bound $C$ exists.

This paper is organized as follows. In Section 2, we provide bounds for the strong odd chromatic number of trees and unicyclic graphs. In Section~3, we consider planar and outerplanar graphs in the same regard and answer Question~\ref{q:constant} in the affirmative. In Section 4, we study strong odd coloring of graph products.

\section{Trees and unicyclic graphs}

Let us call a tree \emph{odd} if the degree of every vertex is odd.

\bpn   \label{t:tree}
For every tree\/ $T$ we have\/ $\soc(T)\leq 3$, and\/ $\soc(T)=2$
 holds if and only if\/ $T$ is odd.
Moreover, a strong odd coloring of $T$ with\/ $\soc(T)$ colors can be
 determined in linear time.
\epn

\begin{proof}
If $T$ is odd, then its unique proper 2-coloring is a strong odd
 coloring, hence $\soc(T)=2$ and of course $T$ can be properly
 2-colored in linear time.
Otherwise two colors are not enough, unless $T$ is trivial, because if a vertex $v$ has a positive even
 degree $\deg_T(v)$, then the closed neighborhood $N_T[v]=N_T(v)\cup\{v\}$ requires at least three colors under any strong odd coloring of $T$.
We give the following linear-time algorithm to show that three
 colors are always enough.

Choose a root vertex $r$ and traverse $T$ with Breadth-First Search.
Let the root $r$ have color 1.
If $\deg_T(r)$ is odd, assign color 2 to all of $N_T(r)$; contrarily,
 color one neighbor of $r$ with 3 and the rest of $N_T(r)$ with 2.
Afterwards, at each non-root vertex $v$---in the order as they are visited
 during Breadth-First Search---if the number of its children is
 even, make the entire $N_T(v)$ monochromatic.
Contrarily, if $v$ has an odd number of children, then let one child
 of $v$ be assigned with the unique color from $\{1,2,3\}$ which is
 distinct from the colors of $v$ and of its parent, and let all the
 other children of $v$ receive the color of the parent of $v$.
At the end of this procedure the entire $T$ is strongly odd $3$-colored.
\end{proof}

\bsk

\bpn
If\/ $G$ is a connected unicyclic  graph other than\/ $C_5$, then
 $\soc(G)\leq 4$. Moreover, a strong odd $4$-coloring of $G$ can be determined in linear time.
\epn

\begin{proof}
It is easy to check  that for cycles we have:

\begin{equation}
    \label{eqn:cycles}
\soc(C_n)=
\begin{cases}
3 & \text{\quad if \,} 3 \mid n \,;\\
4 & \text{\quad if \,} 3 \nmid n \text{ and } n\neq5\,;\\
5 & \text{\quad if \,} n=5 \,.
\end{cases}
\end{equation}
\noindent Moreover, a strong odd $\soc(C_n)$-coloring of $C_n$ can be determined in linear time.

\smallskip

Assume now that $G$ is a connected  unicyclic graph which is not a cycle. Recall that, by Propostion~\ref{t:tree}, every tree can be strongly odd 3-colored in linear time.
Suppose first that the cycle in $G$ is some $C_n$  with $n \neq 5$.   Then, by~\eqref{eqn:cycles}, there exists a strong odd 4-coloring $\varphi$ of the vertices of this cycle. Consider a vertex $a\in C_n$ that has neighbors outside the cycle, and let $b$ and $c$ be the neighbors of $a$ on the cycle. Assume without loss of generality that $\varphi(a)=1$, $\varphi(b)=2$ and $\varphi(c)=3$ (these must be three distinct colors due to the strong coloring of the cycle).
Let $N^*(a)=N_G(a)\backslash\{b,c\}$  be the set of all neighbors of $a$ not on the cycle. If $N^*(a)$ is even-sized, then color it monochromatically with 2. Contrarily, if $N^*(a)$ is odd-sized, color it monochromatically with 4. Note that this partial coloring is still strong odd.  Now each  vertex in $N^*(a)$ is the root of a tree that emanates from it away from the cycle. Each such tree can be colored by three colors avoiding the color 1 (which is used for $a$); so one can use the colors $2,3,4$ and thus complete a strong odd $4$-coloring of $G$. Therefore, unless the cycle of $G$ is $C_5$, four colors suffice for a strong odd coloring of $G$ and such a coloring can be determined in linear time.

Suppose now that the cycle in $G$  is $C_5$ with vertices $a ,b ,c ,d ,e$ in a cyclic order so that the vertex $a$ has neighbors outside the cycle. Let $N^*(a)=N_G(a)\backslash\{b,e\}$. 
Assign colors as follows: $\varphi(a)=1$, $\varphi(b)=2$, $\varphi(c) =3$, $\varphi(d) = 4$, $\varphi(e) = 5$,
the last one considered as a temporary color only.
 Again, we distinguish between two possibilities regarding the parity of $|N^*(a)|$.

\medskip

\noindent \textbf{Case 1:}  \textit{$N^*(a)$ is even-sized}. Split $N^*(a)$ into two odd-sized subsets $X$ and $Y$ (say  $|X| = |N^*(a)| - 1$, $|Y| = 1$). Recolor the vertex $e$ by $2$. Color all of $X$ by $2$ and all of $Y$ by $3$. 
This is still a (partial) strong odd coloring  but with four colors. 
Now apply to  the trees emanating from $N^*(a)$ and from any other vertex of $C_5$ the above argument.

\medskip

\noindent \textbf{Case 2:}  \textit{$N^*(a)$ is odd-sized}. Color all of $N^*(a)$ with $2$. Again, recolor the vertex $e$ by $2$ and then complete a strong odd $4$-coloring as above.

This concludes the proof that apart from $C_5$ every other connected unicyclic graph admits a strong odd $4$-coloring which can be determined in linear time.
\end{proof}

\section{Planar graphs and outerplanar graphs}

\begin{figure}[ht!]
	$$
		\includegraphics[scale=0.5]{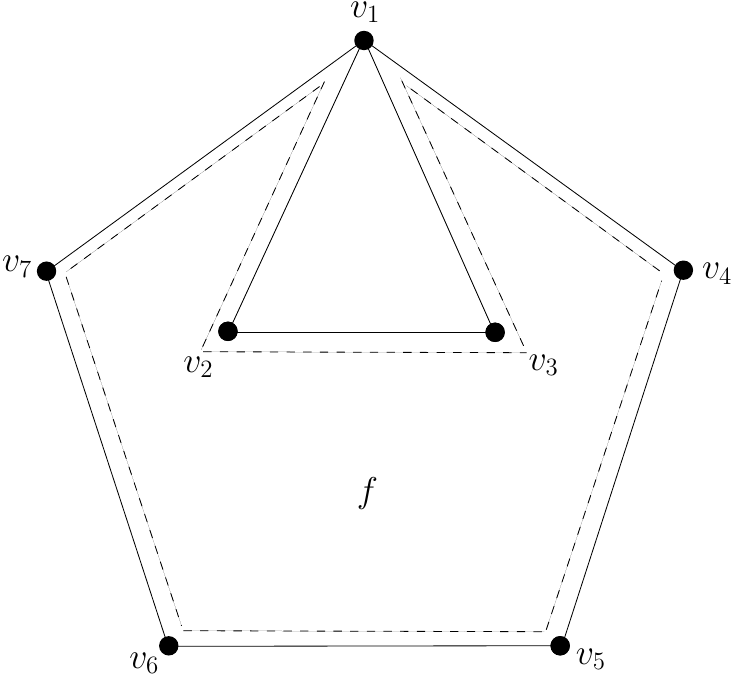}
	$$
	\caption{The boundary walk of the face $f$.}
	\label{fig:faceboundary}
\end{figure}

\def \gta {G^{12.a}}
\def \gtb {G^{12.b}}
Recall that a planar graph is a graph which can be embedded in the (Euclidean) plane. A plane graph is a fixed embedding of a planar graph. 
If a plane graph $G$ is drawn in the plane $\mathcal{P}$, then the maximal connected
regions of $\mathcal{P} \backslash G$ are the \textit{faces} of $G$.  Faces of $G$ are open 2-cells. The \textit{boundary} of a face $f$ is the boundary in the
usual topological sense, denoted $\partial(f)$. It is a collection of all vertices and edges, respectively called \textit{boundary vertices} and \textit{boundary edges} of $f$, contained in the closure of $f$ that can be organized into a
closed walk in $G$ traversing along a simple closed curve lying just inside the face $f$. This closed walk is unique up to the choice
of initial vertex and direction, and is called the \textit{boundary walk} of the face $f$ (see~\cite{GroTuc01}, p.~101).  The sets of boundary vertices and boundary
edges of $f$ are denoted by $V(f)$ and $E(f)$, respectively. Similarly, $F(v)$ denotes the set of faces
incident with a vertex $v$. For all notation defined, the relevant graph may be referred to as a subscript if necessary. When $e_1$ and $e_2$ are a pair of parallel edges which along with their shared endpoints constitute the boundary of a face $f$, such an
$f$ is called a \textit{digon}. Whenever we are dealing with simple graphs, the boundary walk can be seen as a list of consecutive vertices along the boundary of $f$. As an example see the graph in Figure~\ref{fig:faceboundary} where $v_1,v_2,v_3,v_1,v_4,v_5.v_6,v_7,v_1$ is a boundary walk of a face $f$. A graph is \textit{outerplanar} if it can be embedded in the plane such that all vertices are incident with a single face; e.g., the graph in Figure~\ref{fig:faceboundary} is outerplanar.
(A more traditional way of drawing would place the triangle
outside of the pentagon, and an equivalent standard definition of outerplanar graphs requires the existence of a planar embedding where all vertices are on the boundary of the infinite region.)

\begin{figure}[htp!]
    \centering
    \begin{minipage}{0.45\textwidth}
        \centering
        \includegraphics[width=0.9\textwidth]{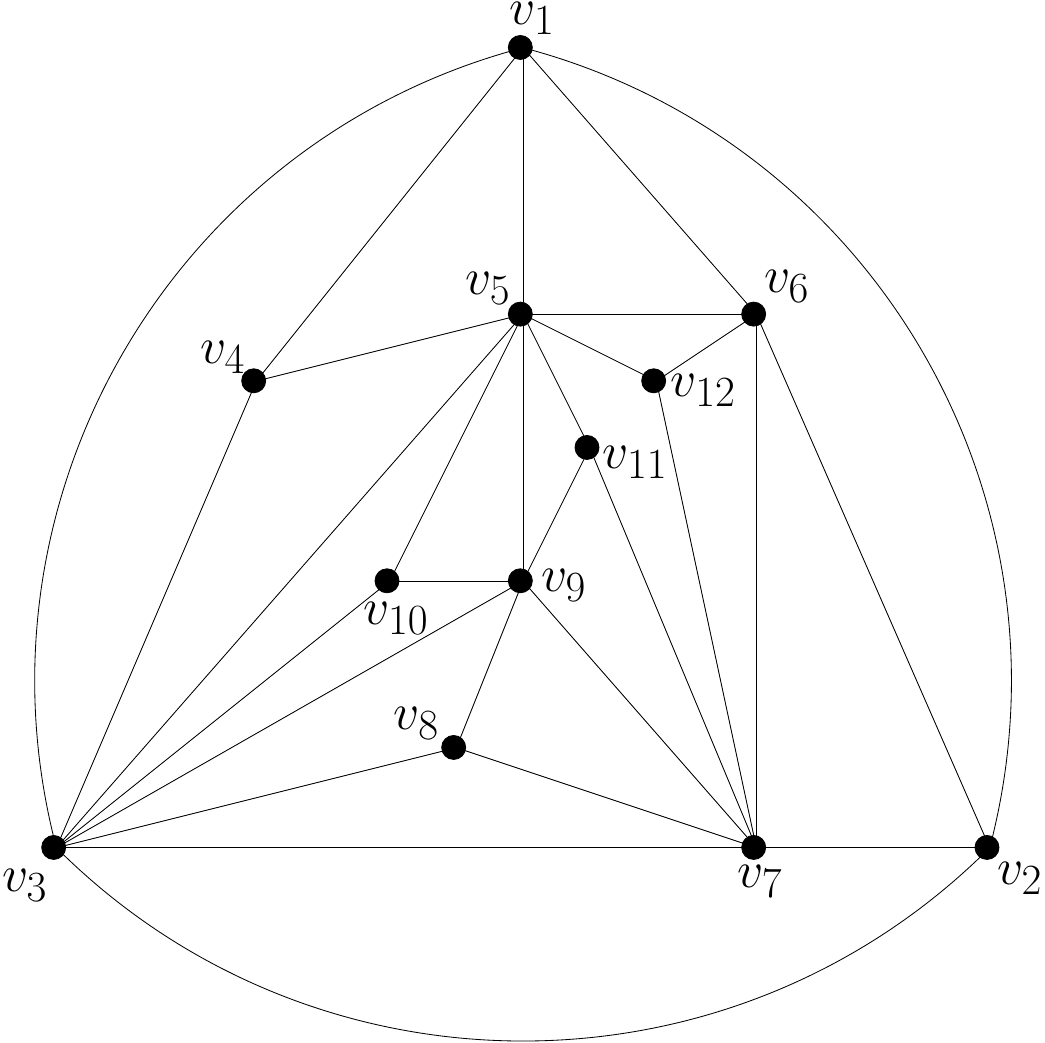}
        \caption{The graph $G^{12.a}$.}
        \label{fig:12a}
    \end{minipage}\hfill
    \begin{minipage}{0.45\textwidth}
        \centering
        \includegraphics[width=0.9\textwidth]{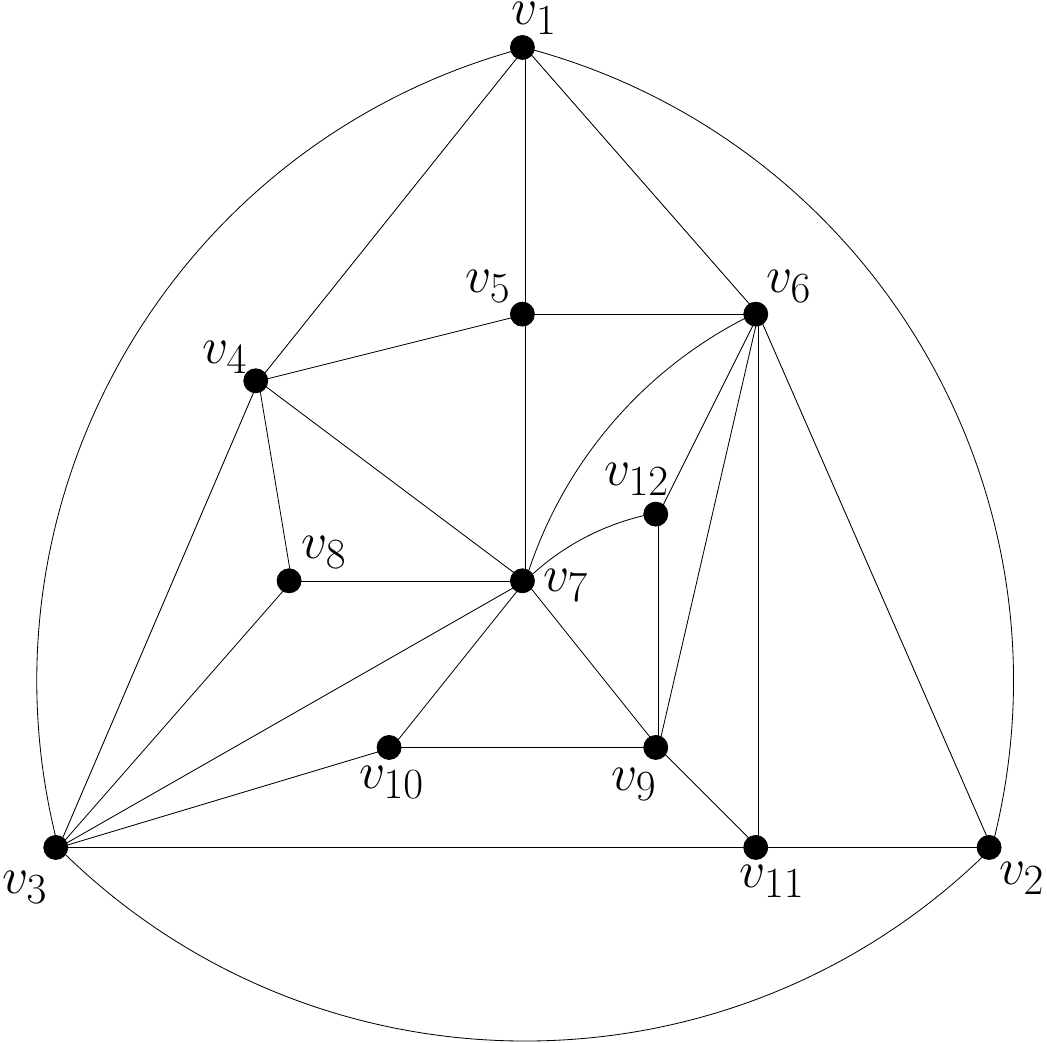}
        \caption{The graph $G^{12.b}$.}
    \end{minipage}
    \label{fig:12b}
\end{figure}

In this section we give estimates on the largest possible value of $\soc$
 in the classes of planar and outerplanar graphs.
Figures~\ref{fig:12a} and~\ref{fig:12b}  exhibit two planar graphs that we denote
 by $\gta$ and $\gtb$, respectively.
Each of those two graphs has order $12$, diameter $2$ and provides the lower bound 12 for
 $\max\lbrace \soc(G) \mid G \mathrm{\ is\ planar} \rbrace$. Later on in this section we shall show that $\max\lbrace \soc(G) \mid G \mathrm{\ is\ planar} \rbrace$ is indeed finite, but first let us discuss $\gta$ and $\gtb$ in regard to their strong odd chromatic number.

\bpn    \label{p:planar}
We have\/ $\soc(\gta) = \soc(\gtb) = 12$.
\epn

\begin{proof}
In both cases the vertices are indexed from 1 to 12.
The color assignment $\vp(v_i)=i$ obviously is a strong odd coloring,
 hence it will suffice to show that no two vertices can be assigned with
 the same color under a strong odd coloring $\varphi$.
We may assume that $\vp(v_i)\leq i$ for all $1\leq i\leq 12$,
 without loss of generality.
(This situation can be attained by renumbering the colors, if necessary. 
Namely, for all colors $i$ simultaneously, for all vertices
 in the color class $\varphi^{-1}(i)$, replace their color 
 with $\min\{j: \varphi(v_j)=i\}$.)
For the duration of this proof, by an ``$i$-neighbor'' we shall mean a neighbor assigned with the color $i$.

\bsk

In the graph $\gta$ we proceed in the natural order $v_1,v_2,\dots,v_{12}$ of vertices, and demonstrate that $\varphi(v_i)=i$ $(1\leq i\leq 12)$.
 \begin{itemize}
  \item $\vp(v_1)=1$, $\vp(v_2)=2$, $\vp(v_3)=3$ --- indeed, the
   three vertices $v_1,v_2,v_3$ induce a $K_3$ in $\gta$.
  \item $\vp(v_4)=4$ --- by adjacency, $\vp(v_4)\notin\lbrace 1,3 \rbrace$;
   $\vp(v_4)\neq 2$, otherwise $v_1$ would need a third 2-neighbor,
   but $N(v_1)=\{v_2,v_3,v_4,v_5,v_6\}$ and $v_5$ is adjacent to $v_4$ whereas $v_6$ is adjacent to $v_2$.
  \item $\vp(v_5)=5$ --- by adjacency, $\vp(v_5)\notin\lbrace 1,3,4 \rbrace$;
   $\vp(v_5)\neq 2$, otherwise $v_1$ would need $\vp(v_6)=2$ but $v_6$ is adjacent to $v_2$ and $v_5$.
  \item $\vp(v_6)=6$ --- by adjacency, $\vp(v_6)\notin\lbrace 1,2,5 \rbrace$;
   $\vp(v_6)\notin\lbrace 3,4 \rbrace$, otherwise $v_1$ would have just two
   3-neighbors or two 4-neighbors.
  \item $\vp(v_7)=7$ --- by adjacency, $\vp(v_7)\notin\lbrace 2,3,6 \rbrace$;
   $\vp(v_7)\notin\lbrace 1,4,5 \rbrace$ because $N(v_3)=\{v_1,v_2,v_4,v_5,v_7,v_8,v_9,v_{10}\})$ and  $v_8,v_9,v_{11},v_{12}$ are
   adjacent to~$v_7$: namely, if $\varphi(v_7)\in\{1,4,5\}$ then $v_3$ would need $\vp(v_{10})=\vp(v_{7})$ to have that color
   a third time in $N(v_3)$, but $(v_5,v_{10})$ is an edge excluding color 5 at $v_{10}$,
   and $\vp(v_{10})\in\lbrace 1,4\rbrace$ would yield only two neighbors of $v_5$ in that color.
  \item $\vp(v_8)=8$ --- by adjacency, $\vp(v_8)\notin\lbrace 3,7 \rbrace$;
   $\vp(v_8)\notin\lbrace 1,2,4 \rbrace$ as otherwise
   $\vp(v_8)=\vp(v_{10})=\vp(v_{11})$ would follow (first for $\vp(v_{10})$
   by vertex $v_3$, second for $\vp(v_{11})$ by vertex $v_9$), but then
   the number of neighbors of $v_5$ and $v_7$ in that color would not have
   the same parity, hence a strong odd coloring could not be completed
   at vertex $v_{12}$;
   $\vp(v_8)=5$ would occur only twice in $N(v_3)$;
   $\vp(v_8)=6$ would imply $\vp(v_{11})=6$ in $N(v_7)$ and then
   $\vp(v_{10})=6$ in both $N(v_5)$ and $N(v_9)$, a contradiction in $N(v_3)$.
  \item $\vp(v_9)=9$ --- by adjacency, $\vp(v_9)\notin\lbrace 3,5,7,8 \rbrace$;
   $\vp(v_9)\notin\lbrace 1,2,4 \rbrace$ due to $N(v_3)$ and the presence
   of the edge $(v_9,v_{10})$;
   $\vp(v_9)\neq 6$ due to $N(v_7)$ as $v_{11}$ is adjacent to $v_9$ and
   $v_{12}$ is adjacent to $v_6$.
  \item $\vp(v_{10})=10$ --- by adjacency, $\vp(v_{10})\notin\lbrace 3,5,9 \rbrace$;
   $\vp(v_{10})\notin\lbrace 1,2,4,7,8 \rbrace$ due to $N(v_3)$;
   $\vp(v_{10})\neq 6$ as otherwise the parity conflict between $N(v_5)$ and $N(v_7)$
   concerning color 6 could not be resolved at vertex $v_{11}$.
  \item $\vp(v_{11})=11$ --- by adjacency, $\vp(v_{11})\notin\lbrace 5,7,9 \rbrace$;
   $\vp(v_{11})\notin\lbrace 1,2,4,8,10 \rbrace$ as any of those colors
   at vertex $v_{11}$ would cause a parity conflict between $N(v_5)$ and $N(v_7)$;
   $\vp(v_{11})\neq 3$ due to $N(v_9)$; $\vp(v_{11})\neq 6$ as otherwise (since $(v_6,v_{12})$ is
   an edge) both $v_5$ and $v_7$ would have just two 6-neighbors.
  \item $\vp(v_{12})=12$ --- every vertex $v_i$ $(1\leq i\leq 11)$ not adjacent to $v_{12}$ has a
   common neighbor $v_j$ with $v_{12}$, hence assigning $\vp(v_{12})=i$ would
   yield that color $i$ occurs precisely twice in $N(v_j)$.
 \end{itemize}

The analysis for $\gtb$ is based on similar principles. Again we proceed in the natural order $v_1,v_2,\dots,v_{12}$ of vertices.
 \begin{itemize}
  \item $\vp(v_1)=1$, $\vp(v_2)=2$, $\vp(v_3)=3$ since $\{v_1,v_2,v_3\}$ induces a $K_3$ in $\gtb$.
  \item $\vp(v_4)=4$ --- by adjacency, $\vp(v_4)\notin\lbrace 1,3 \rbrace$;
   $\vp(v_4)\neq 2$, as otherwise $v_1$ requires the color 2 on either $v_5$ or $v_6$, and both are impossible because of the edges $(v_4,v_5)$ and $(v_2,v_6)$.
  \item $\vp(v_5)=5$ --- by adjacency, $\vp(v_5)\notin\lbrace 1,4 \rbrace$;
   $\vp(v_5)\notin\{2,3\}$ by $N(v_1)$ and the presence of the edge $(v_5,v_6)$.
  \item $\vp(v_6)=6$ --- by adjacency, $\varphi(v_6\notin\{1,2,5\})$; $\varphi\notin\{3,4\}$ by $N(v_1)$.
  \item $\vp(v_7)=7$ --- by adjacency, $\vp(v_7)\notin\lbrace 3,4,5,6 \rbrace$;
   $\vp(v_7)\neq 2$ by $N(v_6)$ and the presence of the edge $(v_2,v_{11})$;
   $\vp(v_7)\neq 1$ by $N(v_5)$.
  \item $\vp(v_8)=8$ --- by adjacency, $\varphi(v_8)\notin\{3,4,7\}$; $\varphi(v_8)\neq6$ as otherwise $\varphi(v_{10})=6$ by $N(v_7)$ and the edges $(v_6,v_{12})$ and $(v_9,v_{12})$, however then $v_9$ would have precisely two $6$-neighbors; $\varphi(v_8)\neq5$ by $N(v_4)$; $\varphi(v_8)\neq2$ as otherwise $\varphi(v_{10})=2$ by $N(v_3)$, but then in view of $N(v_7)$ the color~$2$ is used for either $v_9$ or $v_{12}$, which would force that the color $2$ occurs exactly twice on $N(v_6)$; finally, $\varphi(v_8)\neq1$ by $N(v_4)$.
  \item $\vp(v_9)=9$ --- by adjacency, $\varphi(v_9)\notin\{6,7\}$; $\varphi(v_9)\notin\{4,5,8\}$ by $N(v_7)$ and the edges $(v_9,v_{10})$ and $(v_9,v_{12})$; $\varphi(v_9)\neq3$ by $N(v_{10})$; $\varphi(v_9)\notin\{1,2\}$ by $N(v_6)$ and $(v_9,v_{12})$.
  \item $\vp(v_{10})=10$ --- by adjacency, $\varphi(v_{10})\notin\{3,7,9\}$; $\varphi(v_{10})\neq8$ since otherwise $N(v_3)$ would give $\varphi(v_{11})=8$ and consequently $\varphi(v_{12})\neq 8$ by $N(v_6)$, however that would create for $v_7$ exactly two $8$-neighbors; $\varphi(v_{10})\neq6$ by $N(v_9)$; $\varphi(v_{10})\notin\{4,5\}$ as otherwise $N(v_7)$ would imply that $\varphi(v_{12})=\varphi(v_{10})$, but then that color would appear exactly twice on $N(v_9)$; $\varphi(v_{12})\notin\{1,2\}$ as otherwise $\varphi(v_{11})=\varphi(v_{10})$ by $N(v_3)$, which would give that precisely two neighbors of $v_6$ use the color $\varphi(v_{10})$.
  \item $\vp(v_{11})=11$ --- by adjacency, $\varphi(v_{11})\notin\{2,3,6\}$; $\varphi(v_{11})\notin\{1,4,7,8,10\}$ by $N(v_3)$; $\varphi(v_{11})\notin\{5,9\}$ as otherwise $\varphi(v_{12})=\varphi(v_{11})$ by $N(v_6)$, but then $N(v_7)$ uses the color $\varphi(v_{11})$ exactly twice.
  \item $\vp(v_{12})=12$ --- every vertex $v_i$ $(1\leq i\leq 11)$ not adjacent to $v_{12}$ has a
   common neighbor $v_j$ with $v_{12}$, hence assigning $\vp(v_{12})=i$ would
   yield that color $i$ occurs precisely twice in $N(v_j)$.
 \end{itemize}
\end{proof}

On a similar note, we establish a lower bound for the parameter $\soc$ over the class of outerplanar graphs.

\bpn    \label{p:outerplanar}
The outerplanar graph\/ $G_7=K_1 + P_6$ has\/ $\soc(G_7)=7$.
\epn

\begin{proof}
The vertex from $K_1$, say $v_1$, must have its private color as it is completely adjacent to $P_6$.
If a color appears more than once in $P_6$, then it has to occur
at least three times.
But then there is a vertex in $P_6$ whose both neighbors in $P_6$ have this color,
 and the color of its third neighbor $v_1$ in $G_7$ is distinct, a contradiction.
Thus, each of the seven vertices has its private color.
\end{proof}

Let us now turn to answering the Question~\ref{q:constant}. We provide an affirmative answer by making use of a related concept, initiated in 2009 by Czap and Jendrol’~\cite{CzaJen09}. An equivalent form of their definition reads as follows.
A \textit{facially odd coloring} (resp.\ \textit{proper facially odd coloring}) of a 2-connected plane multigraph $G$  is a coloring (resp.\ proper coloring) of its vertices such that every face is incident with zero or an odd number of vertices of each color. The minimum number of colors in a facially odd (resp.\ proper facially odd) coloring of $G$ is denoted $\chi_{\mathrm{fo}}(G)$ (resp.\ $\chi_{\mathrm{pfo}}(G)$). Clearly, both parameters are well-defined for any plane multigraph; the confinement to $2$-connected multigraphs is due to the original definition of these two coloring notions which instead of incidences of any color to a face counts for the occurrences of any color along each boundary walk of the face. Note that for any plane multigraph $G$ we have $\chi_{\mathrm{fo}}(G)\leq \chi_{\mathrm{pfo}}(G)$.  Czap, Jendro\'{l} and Voigt~\cite{CzaJenVoi11}  proved that every 2-connected plane multigraph $G$ has $\chi_{\mathrm{pfo}}(G)\leq118$. The bound 118 was improved to 97 by Kaiser, Ruck\'{y}, Stehl\'{i}k and \v{S}krekovski~\cite{KaiRucSteSkr14}, which is the current best general upper bound on this graph parameter for the class of $2$-connected plane multigraphs. On the other end of the spectrum, recently, \v{S}torgel~\cite{Sto23} pointed out to an infinite family of $2$-connected plane graphs $G$ with $\chi_{\mathrm{pfo}}(G)=12$; so the general upper bound for $\chi_{\mathrm{pfo}}(G)$ is between 12 and 97. Bounds of $\chi_{\mathrm{pfo}}(G)$ for $3$-connected plane graphs having the property that all faces of a certain size are in a sense far from each other were considered in~\cite{CzaJenKar111}. Wang, Finbow and Wang~\cite{WanFinWan12} proved that only two 2-connected outerplane graphs need 10 colors, whereas the rest admit a proper facially odd coloring with at most 9 colors. As for (possibly improper) facially odd colorings of plane graphs, Czap and Jendrol’~\cite{CzaJen09} conjectured that $\chi_{\mathrm{fo}}(G)\leq 6$ for any
2-connected plane graph $G$, providing an infinite family of examples having $\chi_{\mathrm{fo}}= 6$. 
This conjecture was refuted by Kaiser et al.~\cite{KaiRucSteSkr14}, by showing that the best possible general upper bound on $\chi_{\mathrm{fo}}(G)$
 for the class of 2-connected plane simple graphs $G$ is at least 8.
 We define

 \[\chi_{\mathrm{pfo}}(\mathcal{P}) = \max_{ G\hookrightarrow \mathcal{P}} \chi_{\mathrm{pfo}} (G)\]

\noindent as the maximum of   $\chi_{\mathrm{pfo}}(G)$ over all $2$-connected multigraphs $G$ embedded into the plane $\mathcal{P}$;
and, analogously,
\[
\chi_{\mathrm{pfo}}(\mathcal{P}_\mathrm{o})= \max_{ G\text{ is outerplanar}} \chi_{\mathrm{pfo}}(G)
\]
 for 2-connected outerplanar graphs.
As mentioned in the previous paragraph, it has already  been established that  $12\leq \chi_{\mathrm{pfo}}(\mathcal{P})\leq 97$;
 moreover, the tight result on outerplanar graphs is
  $\chi_{\mathrm{pfo}}(\mathcal{P}_\mathrm{o}) = 10$. 

Note that over the class of 2-connected plane
multigraphs, the best general upper bound for (not necessarily proper) facially odd colorings is again $\chi_{\mathrm{pfo}}(\mathcal{P})$. Indeed, if $G$ is a plane graph that realizes the latter bound in regard to $\chi_{\mathrm{pfo}}$, just replace each edge of $G$ by a digon bounding a face; clearly, any facially odd coloring
of the resulting multigraph is necessarily a proper facially odd coloring of the original graph $G$.

\smallskip

\begin{figure}[ht!]
	$$
		\includegraphics[scale=0.6]{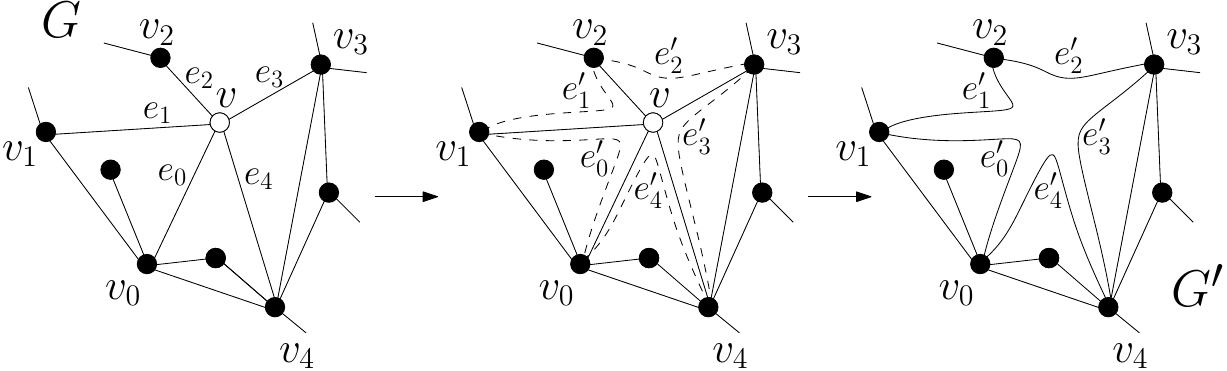}
	$$
	\caption{The annihilation of a vertex $v$. The original graph $G$ is on the left, the
resulting graph $G'$ on the right. The annihilated vertex is depicted as an empty circle.}
	\label{fig:annihilation}
\end{figure}

In the proof of our next result, Theorem~3.4, we shall use the following graph construction (see e.g.~\cite{FabLuzRinSot23, KaiRucSteSkr14}). Given a plane multigraph $G$, let $v \in V_G$ be a vertex of degree $\deg_G(v)=d \geq 2$ which is not incident with a pair of parallel edges. Let the edges incident with $v$ be
enumerated in a clockwise order as $e_i = vv_i$, $i \in \mathbb{Z}_d$ (we write $\mathbb{Z}_d$ for the set $\{0, \ldots, d-1\}$ with addition modulo $d$); thus we assume that the vertices $v_i$ are pairwise different. The \textit{annihilation} of $v$ is the construction of a
plane multigraph $G'$ from $G$ defined as follows:
\begin{enumerate}
\item[(1)] add edges $e'_i = v_iv_{i+1}$, $i\in  \mathbb{Z}_d$, embedded in the plane so that for each $i$, the edges $e_i,e_{i+1}$, and $e'_i$
 constitute a boundary walk of a face;
\item[(2)] remove the vertex $v$ together with all the edges $e_i$.
\end{enumerate}

The assumption that $v$ is not incident with a pair of parallel edges is significant: without
it, the annihilation of $v$ may produce a loop. Intuitively, one may achieve the desired embeddings of the edges $e'_i$ by drawing each $e'_i$ ‘close enough’ to the curve consisting of the embeddings of $e_i$ and $e_{i+1}$; see Figure~\ref{fig:annihilation} for an example of a properly conducted annihilation of a vertex. Notice that if $d=2$ then we introduce two parallel edges, $e'_0$ and $e'_1$, which have as endpoints the pair of vertices $v_0$ and $v_1$.

\smallskip

Regarding the faces of $G$ and $G'$, it is obvious that the following holds:

\bob
    \label{obs}
Let $G'$ be obtained from $G$ by the annihilation of a vertex $v \in V_G$. Then
\begin{enumerate}
\item[(1)] every face of $G$ not in $F_G(v)$ is also a face of $G'$;
\item[(2)] each face $g \in F_G(v)$ has its counterpart $g'$ in $G'$ such that a boundary walk of $g'$ arises from a boundary walk of $g$ by replacing each of its subsequences of the form
$e_ive_{i+1}$ with $e'_i$, and hence $V(g') = V (g)\setminus\{v\}$;
\item[(3)] there is precisely one more face in $G'$, having as
its boundary walk the sequence $v_0e'_
{d-1}v_{d-1}e'_{d-2} \cdots v_1e'_0v_0$.
\end{enumerate}
\eob

We are ready for our main result of this section. It answers in the affirmative the concluding question of Kwon and Park~\cite{KwoPar24}, i.e., Question~\ref{q:constant}.

\btm   \label{t:plane}
\
\begin{itemize}
 \item[$(a)$] 
For every planar graph $G$ it holds that
$$\soc(G)\leq \chi_{\mathrm{pfo}}(\mathcal{P}) \cdot \chi(G)\,.$$
 \item[$(b)$] 
For every outerplanar graph $G$ it holds that
$$\soc(G)\leq \chi_{\mathrm{pfo}}(\mathcal{P}_o) \cdot \chi(G)\,.$$
\end{itemize}
\etm

\begin{proof}
Let $k=\chi(G)$. We consider a plane embedding of $G$, and denote it by $G$ as well.

\bigskip

\noindent \textbf{Claim 1.} \textit{There are plane multigraphs $G_1,G_2,\ldots,G_{k}$ whose vertex sets  partition $V(G)$ into independent subsets such that for any $x\in V_G$ and any $i=1,\ldots,k$ the following holds:
\begin{center}
if $|N_G(v)\cap V(G_i)|\geq2$ then $N_G(v)\cap V_{G_i}=V(f_v)$ for a face $f_v$ of $G_i$.
\end{center}}

\medskip

\noindent \textit{Proof.} Let $\varphi$ be a proper coloring of $G$ with colors $1,2,\ldots,k$ and let $V_i=\varphi^{-1}(i)$ $(i=1,\ldots,k)$ be the corresponding color classes. For any $i=1,\ldots,k$ we construct the multigraph $G_i$ from $G$ as follows:
\begin{itemize}
\item first, remove every edge of $G$ having both endpoints outside $V_i$;
\item next, from the remaining spanning subgraph of $G$ remove any $1^-$-vertex not belonging to $V_i$;
\item finally, in any chosen order, annihilate every vertex not colored by $i$.
\end{itemize}
The resulting plane multigraph is $G_i$ (cf. Figure~\ref{fig:g1}).

\smallskip

In view of Observation~\ref{obs}, the collection $G_1, G_2,\ldots, G_k$ meets all requirements.
\hfill$\diamond$

\begin{figure}[ht!]
	$$
		\includegraphics[scale=0.7]{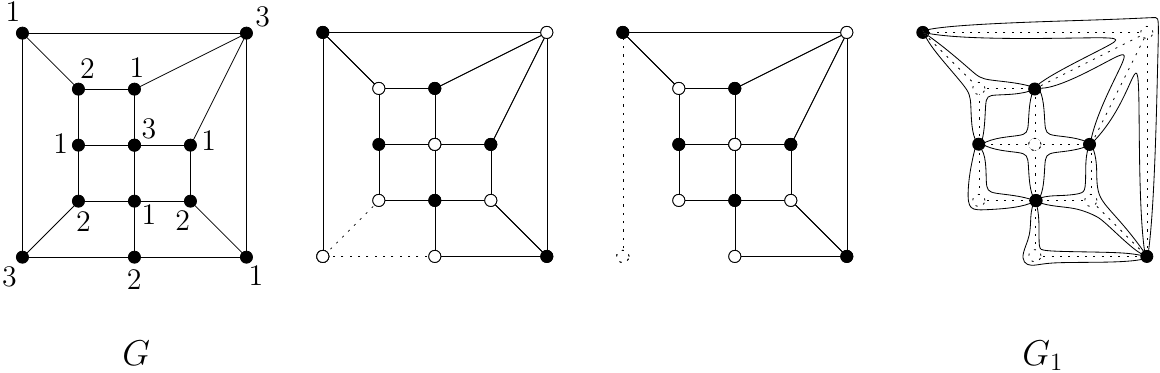}
	$$
	\caption{From left to right, a $3$-chromatic planar graph $G$ with an optimal proper coloring, and the three-step construction of the (multi)graph $G_1$. The empty-circle vertices are the ones not colored by $1$. In each step, the vertices/edges that are being removed appear as dotted.}
	\label{fig:g1}
\end{figure}

\bigskip

\noindent \textbf{Claim 2.} \textit{The edge set of any plane multigraph $H$ on at least three vertices can be enlarged so that the obtained spanning supergraph $\hat{H}$ satisfies the following:
\begin{itemize}
\item[(i)] $\hat{H}$ is a $2$-connected plane multigraph; and
\item[(ii)] For each face $f$ of $H$ there exists a face $\hat{f}$ of $\hat{H}$ with $V(\hat{f})=V(f)$.
\end{itemize}}

\medskip

\noindent \textit{Proof.} First we deal with the possible disconnectedness. Let the input be a plane spanning supergraph $\hat{H}$ of $H$ so that $(ii)$ holds. Initially, we set $\hat{H}=H$. Next, as long as the current input is disconnected, we iterate on the following. There exists a face $f^*$ of $\hat{H}$ whose boundary $\partial(f^*)$ is disconnected. Select a pair of components, say $K'$ and $K''$, of $\hat{H}$ for which it holds that $V_{K'}\cup V_{K''}\subseteq V(f^*)$. Augment $E(\hat{H})$ with a new edge $e$ embedded within the interior of $f^*$ so that $e$ has one endpoint in $K'$ and another endpoint in $K''$; in other words, $e$ is a bridge linking $K'$ and $K''$. Thus $\hat{H}+e$ is a plane spanning supergraph of $H$, and  has one component fewer than $\hat{H}$. Moreover, $\hat{H}+e$ preserves the property $(ii)$; indeed, the addition of the edge $e$ to $E(\hat{H})$ neither destroys an existing face nor creates a new one. The graph $\hat{H}+e$ becomes our new input. Since with each iteration of the described process, the number of components of the input decreases, the end result is a connected plane multigraph, $\hat{H}$, which is a spanning supergraph of $H$ and satisfies $(ii)$.

\begin{figure}[ht!]
	$$
		\includegraphics[scale=0.7]{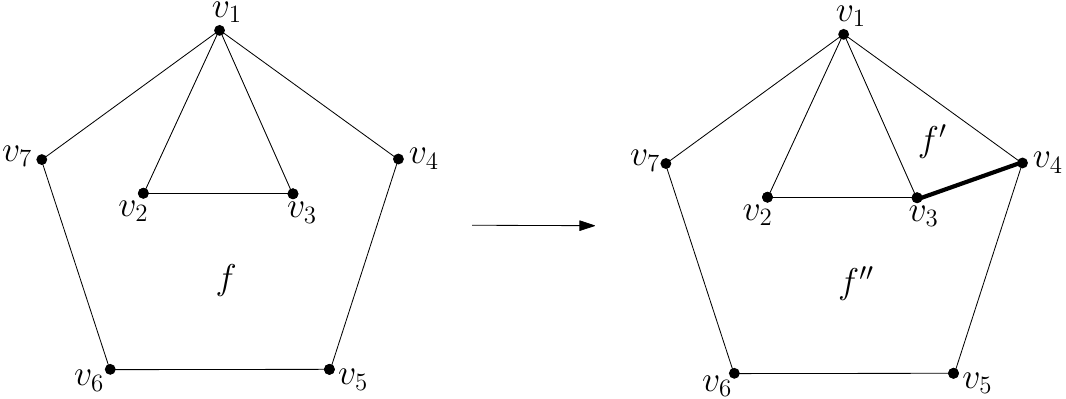}
	$$
	\caption{If $\hat{H}$ is the graph on the left, and $B=v_1v_2v_3$ is the selected end-block (incident with the cut-vertex $v=v_1$) then such a face is $f$. The used boundary walk of $f$ is $v_1v_2v_3v_1v_4v_5v_6v_7v_1$. Here, the vertices $u$ and $w$ are $v_3$ and $v_4$, respectively. We augment $E(\hat{H})$ with the edge $v_3v_4$ (depicted heavier), embedded in the interior of $f$.  The resulting plane graph has one fewer blocks and one more face than $\hat{H}$: namely, the face $f$ splits into a triangular face $f'$ and a face $f''$ such that $V(f'')=V(f)$.}
	\label{fig:endblock}
\end{figure}

Now we take care of the possible cut-vertices of the connected graph $\hat{H}$. As long as there is at least one cut-vertex, select an end-block $B$; say $v$ is the cut-vertex which is incident with $B$. There is a face $f$ of $\hat{H}$ such that $V(B)\cap V(f)\neq\emptyset$ and $V(B)\nsubseteq V(f)$. Thus, along any boundary walk of $f$, the vertex $v$ occurs more than once. There exist vertices $u\in (V(B)\setminus\{v\})\cap V(f)$ and $w\in V(f)\backslash V(B)$ such that a $uv$-edge and a $vw$-edge are facially adjacent in regard to $f$ (i.e., the pair of edges are consecutive on a boundary walk of~$f$).
Augment the edge set of $\hat{H}$ with the (new) edge $uw$ embedded within the interior of $f$.  Thus $\hat{H}+uw$ is a plane graph and has one fewer blocks than $\hat{H}$. Additionally, $\hat{H}+uw$ preserves $(ii)$; indeed, the face sets of $\hat{H}$ and $\hat{H}+uw$ differ only in that the face $f$ of $\hat{H}$ is replaced by two faces $f'$ and $f''$ in $\hat{H}+uw$: namely, $f'$ is a triangular face with $V(f')=\{u,v,w\}$ and $f''$ is a face having the same vertex set as $f$ (cf.\ Figure~\ref{fig:endblock}).  So the property $(ii)$ is preserved in passing from $\hat{H}$ to $\hat{H}+uw$.

We put $\hat{H}+uw$ in place of $\hat{H}$ and iterate. Since with each iteration of the described procedure, the number of blocks of the input decreases, the end result is a $2$-connected plane graph, $\hat{H}$, which is a spanning supergraph of $H$ and satisfies $(ii)$. \hfill$\diamond$

\bigskip

Let $\hat{G_1}, \hat{G_2},\ldots, \hat{G_k}$ be the result of applying Claim~2 to the collection $G_1, G_2,\ldots, G_k$ arising from Claim~1. Take colorings $\varphi_i$ of $\hat{G_i}$ $(i=1,2,\ldots,k)$ with pairwise disjoint color sets so that: if $\hat{G_i}$ is of order at least $3$ then $\varphi_i$ is facially odd, whereas if $\hat{G_i}$ is of order at most $2$ then $\varphi_i$ is proper. The union $\varphi=\varphi_1\cup\varphi_2\cup\cdots\cup\varphi_k$ is a coloring of $G$ using at most $\chi_{\mathrm{pfo}}(\mathcal{P}) \cdot \chi(G)$ colors. Moreover, from Claims 1 and 2 it follows that $\varphi$ is a proper strong odd coloring of $G$.
This completes the proof of $(a)$.

\smallskip

The proof of $(b)$ goes along the same lines, completed by the
 observation that starting from an outerplanar $G$,
 all graphs involved in the argument remain outerplanar.
\end{proof}

So the question raised by the originators of strong odd colorings concerning the existence of a constant bound for the strong odd chromatic number of all planar graphs has an answer in the positive. 
In particular, since $\chi_{\mathrm{pfo}}(\mathcal{P})\leq 97$, the Four Color Theorem of Appel and Haken~\cite{AppHak76,AppHak77,AppHakKoc77} (see also~\cite{AppHak89,RobSanSeyTho93}) and Theorem~\ref{t:plane} combined yield the following corollary.
The second part is derived from Wang, Finbow and Wang's theorem
 $\chi_{\mathrm{pfo}}(\mathcal{P}_{\mathrm{o}})=10$,
 using the fact that every outerplanar graph is 3-colorable.

\bcr    \label{c:plane}
\
\begin{itemize}
 \item[$(a)$] 
For every planar graph $G$ it holds that $\soc(G)\leq 388$.
 \item[$(b)$] 
For every outerplanar graph $G$ it holds that $\soc(G)\leq 30$.
\end{itemize}
\ecr

Let us define
\[
\soc(\mathcal{P}) = \max_{ G\hookrightarrow \mathcal{P}} \soc (G)
\]
as the maximum of  $\soc(G)$ over all planar graphs $G$. In view of Proposition~\ref{p:planar} and Corollary~\ref{c:plane}, we have $12\leq\soc(\mathcal{P})\leq388$. Similarly, let
\[
\soc(\mathcal{P}_\mathrm{o})= \max_{ G\text{ is outerplanar}} \soc (G)
\]
be the maximum of  $\soc(G)$ over all outerplanar graphs $G$.
In view of Proposition~\ref{p:outerplanar} and Corollary~\ref{c:plane}, we have $7\leq \soc(\mathcal{P}_\mathrm{o})\leq 30$.
 We end this section with the following.

 \bpm
 \
 \begin{itemize}
  \item[$(1)$] Is $\soc(\mathcal{P})=12$\,?
  \item[$(2)$] Is $\soc(\mathcal{P}_\mathrm{o})=7$\,?
 \end{itemize}
 \epm

\section{Strong odd colorings of graph products}

In this section we study the behavior of $\soc$ with respect to
 all the four basic types of graph products.

\bdf
Given two graphs $G=(V_G,E_G)$ and $H=(V_H,E_H)$,
 the vertex set in all kinds of their products discussed here
 is the set
  $$
    V_G\times V_H := \{ (g,h) \mid g\in V_G , \, h \in V_H \}
  $$
 (the standard Cartesian product of the corresponding vertex sets).
The edge sets are defined as follows.
 \begin{itemize}
  \item[$(i)$] Cartesian product, $G\Box H$\,:
\begin{eqnarray}
  E(G\Box H) := \{ (g,h)(g',h') & | &
  (g=g'\land hh'\in E_H) \lor (gg'\in E_G \land h=h'), \nonumber \\
  &&
  \ g,g'\in V_G, \ h,h'\in V_H \} .   \nonumber
\end{eqnarray}

  \item[$(ii)$] Direct product, $G\times H$\,:
\begin{eqnarray}
  E(G\times H) :=
  \{ (g,h)(g',h') & | & gg'\in E_G \land hh'\in E_H , \nonumber \\
  &&  g,g'\in V_G, \ h,h'\in V_H \} .   \nonumber
\end{eqnarray}

  \item[$(iii)$] Strong product, $G\boxtimes H$\,:
\begin{eqnarray}
  E(G\boxtimes H) :=
  \{ (g,h)(g',h') & | &
  g'\in N[g] \land h'\in N[h], \nonumber \\
  &&  g,g'\in V_G, \ h,h'\in V_H, \ (g,h) \neq (g',h') \} .   \nonumber
\end{eqnarray}

  \item[$(iv)$] Lexicographic product, $G\circ H$\,:
\begin{eqnarray}
  E(G\circ H) :=
  \{ (g,h)(g',h') & | &
  gg'\in E_G \lor (g=g' \land hh'\in E_H), \nonumber \\
  &&  g,g'\in V_G, \ h,h'\in V_H \} .   \nonumber
\end{eqnarray}

 \end{itemize}
\edf

It turns out that $\soc$ is submultiplicative in all cases,
 with an interesting twist in regard to the $\circ$ operation.

\btm   \label{t:product}
For any two graphs\/ $G,H$ the following inequalities are valid.
 \begin{itemize}
  \item[$(i)$] If\/ $*\in \{ \Box , \times , \boxtimes \}$
   (Cartesian, direct, or strong product), then
   $$\soc(G * H)\leq \soc(G) \, \soc(H).$$

  \item[$(ii)$] In case of the lexicographic product,
   $$\soc(G\circ H) \leq \soc(G) \, (\soc(H+K_1) - 1).$$

 \end{itemize}
\etm

\bpf

As a general notation,  if $F$ is a graph, $v\in V(F)$, and $c$ is
 a color in a vertex coloring of $F$ (any given coloring),
 let us write $d_{c:F}(v)$ for the number of neighbors of $v$
 in color $c$.
If $v=(g,h)$ is a product vertex, we simplify notation and write
 $d_{c:F}(g,h)$ instead of $d_{c:F}((g,h))$.
Similar convention will be applied to color assignments as well,
 writing $\vp(g,h)$ for $\vp((g,h))$.

For any of the four products we let $\vp_G$ be any strong odd
 coloring of $G$ with $\soc(G)$ colors.
Similarly, for the three product types of $(i)$, let $\vp_H$ be any
 strong odd coloring of $H$ with $\soc(H)$ colors.
For $(ii)$ we consider a strong odd coloring of $H+K_1$
 with $\soc(H+K_1)$ colors, and let $\vp_H$ be its restriction
 to $V_H$.
This latter $\vp_H$ uses $\soc(H+K_1)-1$ colors because the color
 of $K_1$ does not appear in $V_H$.

In either case, we compose a coloring of the product graph by the rule
 $$
   \vp(g,h) := ( \vp_G(g) , \vp_H(h) )
 $$
  for all $g\in V_G$, $h\in V_H$.
We claim that every vertex has an odd (or zero) number of neighbors
 in each color.

In $G\Box H$ the neighbors of any vertex $(g,h)$ are of the form
 $(g,h')$ or $(g',h)$, implying that the first or second coordinate of the
 color at each neighbor coincides with the color at the vertex in question.
Moreover, if $(\vp_G(g),\vp_H(h))=(\vp_G(g'),\vp_H(h'))$, then
 $(g,h)$ and $(g',h')$ are not adjacent.
Indeed, if $g=g'$, then $\vp_H(h)=\vp_H(h')$ implies $hh'\notin E_H$;
 and the reasoning is analogous if $h=h'$.
So, for any composite color $c=(c_G,c_H)$,
 the neighbors of any vertex $(g,h)$ assigned with $c$
 belong entirely to a copy of $G$ or to a copy of $H$ in $G\Box H$.
It follows that the number of those neighbors is equal either to
 $d_{c_G:G}(g)$ or to $d_{c_H:H}(h)$, hence odd or zero.

In $G\times H$ if $(g,h)$ and $(g',h')$ are adjacent, we necessarily
 have $\vp_G(g)\neq \vp_G(g')$ and $\vp_H(h)\neq \vp_H(h')$.
Then, for any composite color $c=(c_G,c_H)$, the equality
 $d_{c:G\times H}(g,h) = d_{c_G:G}(g) \, d_{c_H:H}(h)$ is valid.
Consequently $d_{c:G\times H}(g,h)$ is either zero or the product of two odd
 numbers, hence odd.

The open neighborhood $N_{G\boxtimes H}(g,h)$ of a vertex $(g,h)$
 in the strong product is the disjoint union
 $N_{G\Box H}(g,h) \cup N_{G\times H}(g,h)$.
As observed above, the composite color of $(g,h)$ in $G\Box H$
 shares its first or second coordinate with the colors of its
 neighbors, which is not the case with neighbors in $G\times H$.
Thus, $\vp(N_{G\Box H}(g,h)) \cap \vp(N_{G\times H}(g,h))=\empt$,
 and consequently the strong odd coloring property of $\vp$
 is inherited from those of $\vp_G$ and $\vp_H$.
This completes the proof of $(i)$.

For $(ii)$ we first note that
 $$ V(G\circ H) = \bigcup_{g\in V_G} ( \{g\} \times H ) \,;$$
 moreover, if $gg'\in E_G$ then the copies $\{g\} \times H$ and
 $\{g'\} \times H$ of $H$ are completely adjacent, and if
 $gg'\notin E_G$ then they are completely nonadjacent.
Furthermore, in case of adjacency we have
 $\vp(\{g\} \times H) \cap \vp(\{g'\} \times H)=\empt$.

Assume that $(g,h)$ has the composite color $c=(c_G,c_H)$
 in $G\circ H$, and consider its neighbors $(g',h')$ of
 an arbitrary specified color $c'=(c'_G,c'_H)$.
We need to distinguish between two cases.

If $c_G'=c_G$, then we must have $g'=g$.
Indeed, the color class $c_G$ is an independent set in $G$, and
 for $g'\neq g$ the two copies of $H$ are nonadjacent, cannot
 contain any color-$c'$ neighbor of $(g,h)$.
Thus, the number of neighbors of $(g,h)$ in color $c'$ is equal to
 the number of color-$c'_H$ neighbors of vertex $h$ in graph $H$,
  which is odd, due to $\vp_H$.

On the other hand, if $c_G'\neq c_G$, then we claim
 $$
   d_{c':G\circ H}(g,h) = d_{c'_G:G}(g)
     \cdot | \{ h'\in V_H : \vp_H(h') = c'_H \} | \,.
 $$
Indeed, the first term on the right-hand side is the number of
 adjacent copies of $H$ where color $c'_G$ occurs as first coordinate
 of the composite color; this number is odd, due to $\vp_G$.
The second term is of course just the number of vertices in each
 copy of $H$, whose second color-coordinate is $c'_H$.
It is exactly the number of color-$c'_H$ neighbors of $K_1$ in
 $H+K_1$; and this is also odd, because $\vp_H$ is derived from a
 strong odd coloring.
Thus, $d_{c':G\circ H}(g,h)$ is odd, proving $(ii)$.
\epf

\brm
The outerplanar example in Section~3 shows that
 the insertion of a universal vertex may drastically increase
 the strong odd chromatic number; i.e., $\soc(H+K_1)-1$ may be
 substantially larger than $\soc(H)$, and the upper bound
 in $(ii)$ may not be replaced with the formula of $(i)$.
Indeed, $\soc(P_6)=3$, while $\soc(P_6+K_1)=7$.
\erm

Tightness of Theorem \ref{t:product} is clear for strong and
 lexicographic products, since
  $K_p\boxtimes K_q \cong K_p \circ K_q \cong K_{pq}$
 holds for all natural numbers $p$ and $q$.
The situation with the other two types of products is more interesting.
The next theorem disregards the trivial case of $q=1$ since
 $K_p\Box K_1 \cong K_p$ and $K_p\boxtimes K_1$ is the
 edgeless graph $\overline{K}_p \cong p K_1$.

\btm   \label{t:kpkq}
Let\/ $p,q\geq 2$ be integers.
 \begin{itemize}
  \item[$(i)$] We have\/ $\soc(K_p\Box K_q) = pq$ for all\/ $p,q$.
  \item[$(ii)$] If both\/ $p$ and\/ $q$ are odd (and at least\/ $3$),
   then also\/ $\soc(K_p\times K_q) = pq$.
  \item[$(iii)$] If precisely one of\/ $p$ and\/ $q$ is even,
   say\/ $p$ is even and\/ $q$ is odd,
   then\/ $\soc(K_p\times K_q) = q$.
  \item[$(iv)$] If both\/ $p$ and\/ $q$ are even,
   then\/ $\soc(K_p\times K_q) = \min \left( p,q \right)$.
 \end{itemize}
\etm

\bpf
Let $V(K_p) = V_G = \{g_1,\dots,g_p\}$ and
 $V(K_q) = V_H = \{h_1,\dots,h_q\}$.
We refer to the sets $\{g_i\}\times V_H$ ($1\leq i\leq p$) as \emph{columns},
 and to $V_G\times\{h_j\}$ ($1\leq i\leq q$) as \emph{rows}.

In $K_p\Box K_q$ no row and no column can contain any color
 more than once, in any proper vertex coloring.
Moreover, if a color occurs in both column~$i$ and row $j$,
 then vertex $(g_i,h_j)$ has precisely two neighbors of that color
 in $K_p\Box K_q$.
Consequently, no repeated color can occur in any strong odd coloring.
This proves $(i)$.

Turning to $K_p\times K_q$, a vertex coloring is proper if and only if
 each of its color classes is contained entirely in a row or in a column.
We claim the following.
 \begin{itemize}
  \item[$(v)$] \emph{In any strong odd coloring of\/ $K_p\times K_q$,
   if\/ $p$ is odd, then no row can contain a color more than once;
   and similarly, if\/ $q$ is odd, then no column can contain a color more than once.}
 \end{itemize}

Since $K_p\times K_q \cong K_q\times K_p$, it suffices to prove
 $(v)$ for rows.
Suppose for argument's sake that for some $r\geq 2$ the vertices
 $(g_1,h_1),(g_2,h_1),\dots,(g_r,h_1)$ have the same color $c$, but
 none of the vertices $(g_\ell,h_1)$ have color $c$ for $r<\ell\leq p$.
If $r=p$, then vertex $(g_1,h_2)$ has exactly $p-1$ neighbors
 of color $c$, which is even.
Otherwise, if $r<p$, then $(g_r,h_2)$ has $r-1\geq 1$ neighbors of
 color $c$, and  $(g_{r+1},h_2)$ has $r$ neighbors of color $c$;
 and one of them is even.
Hence, in either case color $c$ occurs an even number of times in the
 neighborhood of at least one vertex, proving $(v)$.

Observe that $(v)$ implies $(ii)$, because if both $p$ and $q$ are
 odd, then no color can be repeated at all.

If $q$ is odd, then $(v)$ implies that exactly $q$ colors are
 used in each column.
But for $p$ even, $q$ colors are also enough, because we can then
 make each row monochromatic, assigning color $j$ to all vertices
 in $V_G\times\{h_j\}$ ($j=1,2,\dots,q$).
In that case each vertex outside of row $j$ has exactly $p-1$
 neighbors of color~$j$.
This proves $(iii)$.

Finally, if both $p$ and $q$ are even, then making each row or
 each column monochromatic, we obtain a strong odd coloring with
 $\min\left(p,q\right)$ colors.
This is optimal because $\soc(K_p\times K_q) \geq
 \chi(K_p\times K_q) \geq pq / \alpha(K_p\times K_q) =
 pq / \max\left(p,q\right) = \min\left(p,q\right)$.
This proves $(iv)$.
\epf

\brm
The parameter $\soc$ is not multiplicative with respect to the
 Cartesian product operation.
As an example, consider $F=G\Box H$, where
 $G=g_1g_2g_3g_4g_5\cong C_5$ and $H=h_1h_2h_3h_4h_5\cong C_5$.
The following table exhibits a strong odd coloring of $F$
 with five colors.

\renewcommand{\arraystretch}{1.6}
\begin{center}
 \begin{tabular}{c||c|c|c|c|c|}
     &$h_1$&$h_2$&$h_3$&$h_4$&$h_5$ \\
        \hline
        \hline
   $g_1$  &1&2&3&4&5 \\
        \hline
   $g_2$  &4&5&1&2&3 \\
        \hline
   $g_3$  &2&3&4&5&1 \\
        \hline
   $g_4$  &5&1&2&3&4 \\
        \hline
   $g_5$  &3&4&5&1&2 \\
        \hline
 \end{tabular}
\end{center}
Each cell $(g_i,h_j)$ indicates the color of the corresponding vertex.
Here the closed neighborhood of each vertex contains each color exactly once. \hfill $\Box$
\erm

\subsection{Nordhaus--Gaddum-type consequences}

From the above we can derive that there do not exist any
 nontrivial lower or upper bounds for the strong odd chromatic number of a graph and its complement
 as a function of the number of vertices.
We express this fact in the following result.

\btm
Let\/ $n=(2k+1)^2$ be an odd square.
 \begin{itemize}
  \item[$(i)$] There exists a graph\/ $H_1$ of order\/ $n$ such that
   both\/ $\soc(H_1)=\sqrt n$ and\/ $\soc(\overline{H_1})=\sqrt n$.
  \item[$(ii)$] There exists a graph\/ $H_2$ of order\/ $n$ such that
   both\/ $\soc(H_2)=n$ and\/ $\soc(\overline{H_2})=n$.
  \item[$(iii)$] For every graph\/ $G$ on\/ $n$ vertices,
    $$
      2 \sqrt n \leq \soc(G) + \soc(\overline{G}) \leq 2n \,,
    $$
    $$
      n \leq \soc(G) \cdot \soc(\overline{G}) \leq n^2 ,
    $$
   and all these inequalities are tight whenever\/ $n$ is an odd square.
 \end{itemize}
\etm

\bpf
We can take $H_1 = (2k+1) K_{2k+1}$, whose complement is the
 complete $(2k+1)$-partite graph $K_{2k+1,\dots,2k+1}$.
Each component $K_{2k+1}$ of $H_1$ has a strong odd coloring from
 the color set $\{1,2,\dots,2k+1\}$.
Also, assigning color $i$ to all vertices of the $i^{\mathrm{th}}$
 component, a strong odd coloring of $\overline{H_1}$ with
 $2k+1 = \sqrt n$ colors is obtained.
This proves $(i)$.

Further, we can take $H_2 = K_{2k+1} \Box K_{2k+1}$, whose complement
 is $\overline{H_2} = K_{2k+1} \times K_{2k+1}$.
Then, by parts $(i)$ and $(ii)$ of Theorem \ref{t:kpkq}, we have
 $\soc(H_2) = \soc(\overline{H_2}) = n$.
This proves $(ii)$.

The lower bounds in $(iii)$ follow from the additive and multiplicative
 versions of the Nordhaus--Gaddum theorem on the chromatic number,
 since $\soc(G) \geq \chi(G)$ holds
 for all graphs $G$ by definition.
The upper bounds in $(iii)$ are clear as the obvious
 inequality $\soc(G)\leq |V(G)|$ is valid for all $G$.
Tightness follows from the constructions given for $(i)$ and $(ii)$.
\epf

\brm
Sometimes $\soc(G)$ behaves  nearly as $\chi(G^2)$, as shown
 by the following two examples.
 \begin{itemize}
  \item The product graph $G=K_p\Box K_p$ with $n=p^2$ vertices
   satisfies $\soc(G) + \soc(\overline{G}) = 2n$.
  We also note that it is regular of degree $2p-2$,
   hence $\soc(G) = p^2  \approx (\Delta(G))^2/4$.
  \item Strongly regular graphs $G\in \mathrm{Srg}(n, k , \lambda , \mu )$
   with appropriate choices of the parameters have the property that
   both $G$ and $\overline{G}$ have diameter $2$, hence if $|G|=n$
   then $\chi(G^2) + \chi((\overline{G})^2) = 2n$.
 Note, however, that $G^2$ and $(\overline{G})^2$ are
 not complementary to each other.
 \end{itemize}
\erm

We close this paper with the following question.

\bpm
Does $\soc(G) \leq (1+o(1))\Delta^2\!/4$ hold for every
 graph $G$ of maximum degree $\Delta$, as $\Delta$ gets large?
\epm

\paragraph{Acknowledgements.}

This research of the fourth author was supported in part
 by the National Research, Development and Innovation Office,
 NKFIH Grant FK 132060.

\end{document}